\documentclass{amsart}

\usepackage{graphics}

\newcommand{\mathsym}[1]{{}}

\usepackage{amssymb,amscd}
\usepackage{amsmath}
\usepackage{amsthm}

\newtheorem{theorem}{Theorem}
\newtheorem{lemma}[theorem]{Lemma}

\newtheorem{proposition}[theorem]{Proposition}

\numberwithin{theorem}{section}

\theoremstyle{definition}

\newtheorem{example}[theorem]{Example}

\theoremstyle{remark}

\numberwithin{equation}{section}

\newfont{\germ}{eufm10}

\newcommand{\Z}{{\mathbb Z}}
\newcommand{\R}{{\mathbb R}}

\begin{document}

\begin{center}
{\Large{\bf The Bethe ansatz in a periodic box-ball system}}
\vspace{1mm}\\
{\Large{\bf and }}
\vspace{1mm}\\
{\Large{\bf the ultradiscrete Riemann theta function}}
\vspace{6mm}\\
{\large Atsuo Kuniba and Reiho Sakamoto}
\end{center}

\vspace{0.2cm}
\begin{quotation}
{\small
A}{\tiny BSTRACT.}
{\small
Vertex models with quantum group symmetry 
give rise to integrable cellular automata at $q=0$.
We study a prototype example known as 
the periodic box-ball system.
The initial value problem is solved 
in terms of an ultradiscrete analogue of the
Riemann theta function whose period matrix 
originates in the Bethe ansatz at $q=0$.}
\end{quotation}

\section{Introduction}\label{sec:intro}
The periodic box-ball system \cite{KTT,YYT} is a
completely integrable one-dimensional cellular automaton.
Its dynamics is described as a motion of balls hopping   
exclusively along the periodical 
array of boxes having capacity 1.
The system is identified with a
solvable vertex model \cite{Ba} associated with 
quantum affine algebra $U_q(\widehat{sl}_2)$ at $q=0$, 
where the fusion transfer matrices $T_1, T_2, \ldots$ yield a 
commuting family of deterministic time evolutions.

In \cite{KTT}, the initial value problem 
of the periodic box-ball system is solved by an 
inverse scattering method.
It is done by synthesizing 
the combinatorial versions of the Bethe ansatz \cite{Be} 
at $q=1$ \cite{KR} and $q=0$ \cite{KN}.
The action-angle variables are introduced 
by generalizing the
rigged configurations ($q=1$) up to some 
equivalence specified by the string center equation 
($q=0$).
It enables one to determine the 
time evolution $T^t_l(p)$ of any state $p$
by an explicit algorithm 
whose computational steps are independent of 
the time $t$.

The Bethe ansatz approach \cite{KTT} captures several 
characteristic features 
in the quasi-periodic solutions of soliton equations \cite{DT,DMN}.
For instance, the original nonlinear dynamics 
becomes a straight motion of the Bethe roots (angle variable)
which live in an ultradiscrete analogue 
(\ref{eq:udj}) of the Jacobi variety.

In this paper we exploit such an analogy further 
by representing the solution of the 
initial value problem explicitly in terms of the ultradiscretization (UD) 
of the Riemann theta function (${\bf z} \in \R^g$):
\begin{equation}\label{eq:urt}
\begin{split}
\Theta({\bf z}) &= \lim_{\epsilon \rightarrow +0}
\epsilon\log\left( \sum_{{\bf n} \in \Z^g}
\exp\Bigl(-\frac{{}^t{\bf n}A{\bf n}/2+{}^t{\bf n}{\bf z}}{\epsilon}
\Bigr)\right)\\
&= -\min_{{\bf n} \in \Z^g}
\{{}^t{\bf n}A{\bf n}/2+{}^t{\bf n}{\bf z}\}.
\end{split}
\end{equation}
Here $A$
is the symmetric positive definite $g \times g$ integer matrix
(\ref{eq:A}) 
appearing in the string center equation (\ref{eq:sce}) 
introduced in \cite{KN}.
Likewise the Riemann theta function, 
$\Theta({\bf z})$ enjoys the quasi-periodicity: 
\begin{equation}\label{eq:qp}
\Theta({\bf z} + {\bf v}) = 
{}^t{\bf v}A^{-1}({\bf z} + {\bf v}/2) + \Theta({\bf z})\quad 
\hbox{for any } \; {\bf v} \in 
\Gamma = A\Z^g.
\end{equation}

Let $c_L({\bf n})={}^t{\bf n}A{\bf n}/2 + {}^t{\bf n}{\bf z}$
be the quadratic form appearing in (\ref{eq:urt}),
where $L$ denotes the system size that enters $A$ 
and ${\bf z}$ in our main formula (\ref{eq:y}).
The ultradiscrete Riemann theta function $\Theta({\bf z})$ 
can be spotted in the following degeneration scheme: 
\begin{equation}\label{eq:4}
\begin{split}
& \qquad\qquad\qquad \qquad
\sum_{{\bf n} \in \Z^g}
\exp\bigr(-c_L({\bf n})/\epsilon\bigl)\\
& \quad\qquad \vspace{-0.1cm}{L \rightarrow \infty} \swarrow
\qquad \qquad\qquad\qquad \qquad\searrow {\rm UD} \\
&\sum_{{\bf n} \in \{0,1\}^g}
\exp\bigr(-c({\bf n})/\epsilon\bigl)
\qquad\qquad\qquad\qquad
-\min_{{\bf n} \in \Z^g}\{c_L({\bf n})\} = \Theta({\bf z})\\
&  \quad \qquad\qquad \vspace{-0.1cm}{\rm UD} \searrow
\qquad \qquad\qquad\qquad \qquad\swarrow 
L \rightarrow \infty\\
&\qquad\qquad\qquad\qquad\quad
-\min_{{\bf n} \in \{0,1\}^g}\{c({\bf n})\}
\end{split}
\end{equation}
At the top there is the Riemann theta function, 
which degenerates into various objects.
The UD procedure (\ref{eq:urt}) for getting 
$\Theta({\bf z})$ is 
the SE arrow from the top. 
Then in the limit $L \rightarrow \infty$, the minimum over 
${\bf n} \in \Z^g$ shrinks down to that over 
${\bf n} \in \{0,1\}^g$, which reduces $c_L({\bf n})$ to 
its $L$-independent part $c({\bf n})$.
Consequently, $\Theta({\bf z})$ tends to 
the bottom one in (\ref{eq:4}),
which we call the {\em ultradiscrete tau function}.
The resulting expression (\ref{eq:tau}) for the infinite system 
gives the piecewise linear formula for 
the Kerov-Kirillov-Reshetikhin (KKR) bijection \cite{KR}
from rigged configurations to highest paths. 
One may go down the diagram (\ref{eq:4}) 
via the other route.
The thereby encountered function in the middle left is the 
sum of $2^g$ ``trigonometric terms" that are 
characteristic in the tau functions of soliton solutions 
for the infinite system \cite{JM}.
In fact a procedure analogous to the SW arrow from the top 
has been described in p3.253 in \cite{M}, where quasi-periodic 
soliton solutions tend to those in the infinite system.  

In our approach,  the ultradiscrete 
Riemann theta function ${\Theta}({\bf z})$ arises most naturally
by going from the bottom in (\ref{eq:4}) into the NE direction.
The essential idea \cite{KTT} is to embed a state $p$ of the 
periodic box-ball system into an infinite system as 
$p \otimes p \otimes p \otimes \cdots$.
It turns out that 
the ultradiscrete tau function for such periodic states
is nothing but ${\Theta}({\bf z})$ 
up to irrelevant contributions.
As an application we extend the problem to 
$({\mathbb C}^2)^{\otimes L}$ and 
construct joint eigenvectors
of the commuting time evolutions.
The result may be viewed as an explicit formula 
of the Bethe vectors at $q=0$ 
in terms of the ultradiscrete Riemann theta function. 

In section \ref{sec:pbbs}, we recall the periodic box-ball system 
and the inverse scattering algorithm 
that solves the initial value problem \cite{KTT}.
Section \ref{sec:main} contains our main theorem \ref{th:main}.
Section \ref{sec:discussion} gives the discussion on the 
connection with the Bethe ansatz at $q=0$ \cite{KN}.

We did not intend to make the paper completely self-contained.
Exposition of the KKR bijection \cite{KR} 
and Lemma \ref{lem:waru} have been attributed to \cite{KTT}.
Rather, we have employed a casual description 
to clarify how the algorithmic solution 
to the initial value problem \cite{KTT} 
leads directly to the explicit formula (\ref{eq:y}).  
We shall exclusively consider the case 
where the amplitudes of the solitons are all distinct,
which greatly simplifies the presentation.  
The general case can be treated with the same idea.

\section{Periodic box-ball system and inverse scattering transform}
\label{sec:pbbs}
Let us quickly recall the periodic box-ball system without getting 
much into the crystal base theory. For a comprehensive 
treatment, see \cite{KTT}.
For a positive integer $l$, let 
$B_l =\{(x_1,x_2) \in (\Z_{\ge 0})^2 \mid x_1+x_2=l\}$ and 
set $u_l = (l,0) \in B_l$. 
The two elements $(1,0)$ and $(0,1)$ in $B_1$ will be denoted by
$1$ and $2$ for short. (Thus $u_1=1$.)
In the following, the symbol $\otimes$ meaning the tensor product 
of crystals can just be understood as a product of sets.
Define the map $R: B_l \otimes B_1 \rightarrow B_1 \otimes B_l$ by
\begin{align*}
(x_1,x_2)\otimes 1 &\mapsto 
\begin{cases}1 \otimes (l,0) &\hbox{if } (x_1,x_2)=(l,0)\\
2 \otimes (x_1+1,x_2-1) & \hbox{otherwise},
\end{cases}\\
(x_1,x_2)\otimes 2 &\mapsto 
\begin{cases}2 \otimes (0,l) &\hbox{if } (x_1,x_2)=(0,l)\\
1 \otimes (x_1-1,x_2+1) & \hbox{otherwise}.
\end{cases}
\end{align*}
$R$ is a bijection and called the combinatorial $R$.
We write the relation $R(u \otimes b)= b' \otimes u'$ 
simply as $u \otimes b \simeq b' \otimes u'$, and similarly for 
any consequent relation of the form
$a\otimes u \otimes b \otimes c \simeq
a \otimes b' \otimes u' \otimes c$.

A state of the periodic box-ball system 
is an array of $1$ and $2$, which is 
regarded as an element 
$b_1 \otimes \cdots \otimes b_L \in B_1^{\otimes L}$
with $L$ being the system size.
Let the number of $2 \in B_1$ appearing in 
$b_1 \otimes \cdots \otimes b_L$ be $M$.
Without loss of generality we assume $L \ge 2M$
(see \cite{KTT}, section 3.3).
Let ${\mathcal P}$ be the set of such states.
Then the time evolution 
$T_l: {\mathcal P} \rightarrow {\mathcal P}$
is defined by
\begin{equation}\label{eq:tl}
u_l \otimes p \simeq p^\ast \otimes v_l,\quad 
v_l \otimes p \simeq T_l(p)\otimes v_l.
\end{equation}
In the first relation, one applies the combinatorial $R$ for 
$L$ times to carry $u_l$ through $p \in {\mathcal P}$ 
to the right. 
This determines $v_l \in B_l$ and 
$p^\ast \in {\mathcal P}$ uniquely. 
($p^\ast$ does not play an essential role.)
Then the second relation using the so obtained $v_l$ 
specifies $T_l(p)$, where 
the appearance of the same $v_l$ in the right hand side is a
non-trivial claim (\cite{KTT}, section 2.2).
$v_l$ is dependent on $p$ as opposed to $u_l$.

The combinatorial $R$ is the identity map on $B_1 \otimes B_1$, 
and therefore
$T_1$ is just the cyclic shift 
$T_1(b_1 \otimes \cdots \otimes b_L) = 
b_L \otimes b_1 \otimes \cdots \otimes b_{L-1}$.
The commutativity $T_lT_k = T_kT_l$ holds for any $k, l$
(\cite{KTT}, Theorem 2.2).
\begin{example}\label{ex:t23}
The time evolutions $p, T_l(p), \ldots, T^9_l(p)$ 
of the state $p$ on the top line are listed downward
for $l=2$ and $3$. The system size is  $L=14$.
We omit the symbol $\otimes$.
\[
\begin{array}{llll}
\quad \;\hbox{evolution under } T_2 
\quad \qquad \quad\qquad \hbox{evolution under } T_3 \\
 1 \; 1 \; 2 \; 1 \; 1 \; 1 \; 2 \; 2 \; 2 \; 1 \; 1 \; 1 \; 2 \; 2 \; 
\qquad\quad \; 1 \; 1 \; 2 \; 1 \; 1 \; 1 \; 2 \; 2 \; 2 \; 1 \; 1 \; 1 \; 2 \; 2 \\
 2 \; 2 \; 1 \; 2 \; 1 \; 1 \; 1 \; 1 \; 2 \; 2 \; 2 \; 1 \; 1 \; 1 \; 
\qquad\quad \; 2 \; 2 \; 1 \; 2 \; 1 \; 1 \; 1 \; 1 \; 1 \; 2 \; 2 \; 2 \; 1 \; 1 \\
 1 \; 1 \; 2 \; 1 \; 2 \; 2 \; 1 \; 1 \; 1 \; 1 \; 2 \; 2 \; 2 \; 1 \; 
\qquad\quad \; 1 \; 1 \; 2 \; 1 \; 2 \; 2 \; 2 \; 1 \; 1 \; 1 \; 1 \; 1 \; 2 \; 2 \\
 2 \; 1 \; 1 \; 2 \; 1 \; 1 \; 2 \; 2 \; 1 \; 1 \; 1 \; 1 \; 2 \; 2 \; 
\qquad\quad \; 2 \; 2 \; 1 \; 2 \; 1 \; 1 \; 1 \; 2 \; 2 \; 2 \; 1 \; 1 \; 1 \; 1 \\
 2 \; 2 \; 2 \; 1 \; 2 \; 1 \; 1 \; 1 \; 2 \; 2 \; 1 \; 1 \; 1 \; 1 \; 
\qquad\quad \; 1 \; 1 \; 2 \; 1 \; 2 \; 2 \; 1 \; 1 \; 1 \; 1 \; 2 \; 2 \; 2 \; 1 \\
 1 \; 1 \; 2 \; 2 \; 1 \; 2 \; 2 \; 1 \; 1 \; 1 \; 2 \; 2 \; 1 \; 1 \; 
\qquad\quad \; 2 \; 2 \; 1 \; 2 \; 1 \; 1 \; 2 \; 2 \; 1 \; 1 \; 1 \; 1 \; 1 \; 2 \\
 1 \; 1 \; 1 \; 1 \; 2 \; 1 \; 2 \; 2 \; 2 \; 1 \; 1 \; 1 \; 2 \; 2 \; 
\qquad\quad \; 1 \; 1 \; 2 \; 1 \; 2 \; 2 \; 1 \; 1 \; 2 \; 2 \; 2 \; 1 \; 1 \; 1 \\
 2 \; 2 \; 1 \; 1 \; 1 \; 2 \; 1 \; 1 \; 2 \; 2 \; 2 \; 1 \; 1 \; 1 \; 
\qquad\quad \; 1 \; 1 \; 1 \; 2 \; 1 \; 1 \; 2 \; 2 \; 1 \; 1 \; 1 \; 2 \; 2 \; 2 \\
 1 \; 1 \; 2 \; 2 \; 1 \; 1 \; 2 \; 1 \; 1 \; 1 \; 2 \; 2 \; 2 \; 1 \; 
\qquad\quad \; 2 \; 2 \; 2 \; 1 \; 2 \; 1 \; 1 \; 1 \; 2 \; 2 \; 1 \; 1 \; 1 \; 1 \\
 2 \; 1 \; 1 \; 1 \; 2 \; 2 \; 1 \; 2 \; 1 \; 1 \; 1 \; 1 \; 2 \; 2 \; 
\qquad\quad \; 1 \; 1 \; 1 \; 2 \; 1 \; 2 \; 2 \; 2 \; 1 \; 1 \; 2 \; 2 \; 1 \; 1
\end{array}
\]
Regarding $1$ as an empty box and $2$ as a ball, 
these patterns exhibit the nonlinear dynamics of balls.
There are three solitons (wavepackets) 
with amplitudes $3,2$ and $1$ 
traveling to the right.
\end{example}

Let us proceed to the direct and inverse scattering transforms.
A state $p=b_1\otimes \cdots \otimes b_L$ is called {\em highest} if
\begin{equation*}
\sharp\{1 \le i \le k \mid b_i=1\} \ge 
\sharp\{1 \le i \le k \mid b_i=2\}
\quad \hbox{for all }\; 1 \le k \le L.
\end{equation*}
The state on the top line in example \ref{ex:t23} is highest,
whereas those on the second lines are not.
Let ${\mathcal P}_+$ be the subset of ${\mathcal P}$
consisting of the highest states.
Any state $p \in {\mathcal P}$ can be expressed 
as $p = T_1^d(p_+)$ using some $d \in \Z$ and a highest state 
$p_+ \in {\mathcal P}_+$. 
For instance, the state $T_2(p)$ in example \ref{ex:t23} is written as
$22121111222111=T_1^2(12111122211122)$.
Given a state $p$, such $d$ and $p_+$ are not unique in general.
Picking any one of them will be denoted by $p \mapsto (d, p_+)$.
Consider the KKR bijection $\phi$ from the highest state $p_+$ to 
the rigged configuration \cite{KR}:
\begin{equation}\label{eq:phi}
\begin{picture}(60,53)(-40,-50)
\put(0,0){\line(0,-1){40}}
\put(40,-30){\line(0,-1){10}}\put(42.5,-38.5){$\scriptstyle{J_{i_1}}$}
\put(50,-20){\line(0,-1){10}}\put(52.5,-28.5){$\scriptstyle{J_{i_2}}$}
\put(70,-10){\line(0,-1){10}}
\put(80,0){\line(0,-1){10}}\put(82.5,-8.5){$\scriptstyle{J_{i_g}}$}

\put(-30,-18){$\scriptstyle{\phi}$}
\put(-55,-25){$p_+ \;\;\longmapsto$}

\put(0,0){\line(1,0){80}}
\put(14,-5){\vector(-1,0){13}}
\put(18,-7.5){$\scriptstyle{i_g}$}
\put(26,-5){\vector(1,0){53}}

\put(70,-10){\line(1,0){10}}
\put(18,-17.5){$\scriptstyle{\cdots}$}

\put(50,-20){\line(1,0){20}}

\put(40,-30){\line(1,0){10}}
\put(14,-25){\vector(-1,0){13}}
\put(18,-27.5){$\scriptstyle{i_2}$}
\put(26,-25){\vector(1,0){23}}

\put(0,-40){\line(1,0){40}}
\put(14,-35){\vector(-1,0){13}}
\put(18,-37.5){$\scriptstyle{i_1}$}
\put(26,-35){\vector(1,0){13}}

\end{picture}
\end{equation} 
 
\noindent
The partition $(i_g, \ldots, i_2, i_1)$ is called 
the configuration and the integers $0 \le J_i \le p_i$ are called
the rigging. The combined data define a rigged configuration.
Here $p_i$ is the vacancy number:
\begin{equation}\label{eq:pi}
p_i= L - 2\sum_{j \in \mu}\min(i,j),
\end{equation}
where $\mu = \{i_1 \le i_2 \le  \cdots \le  i_g\}$.
Obviously, $p_{i_1}\ge p_{i_2} \ge \cdots \ge p_{i_g}$ holds,
and it is known that
$i_1 + \cdots + i_g$ coincides with the
number $M$ of $b_k=2 \in B_1$ contained in
$p_+ = b_1\otimes \cdots \otimes b_L$.
Thus we have $p_{i_g} = L-2M \ge 0$ by the assumption. 
See appendix A in \cite{KTT} for an exposition 
adapted to the present context.

The configuration $\mu$ is 
actually independent of the non-uniqueness of the 
choice of $p_+$, and determined solely from $p$. 
The states are classified according to their configurations:
\begin{equation*}
{\mathcal P} = \bigsqcup_{\mu}{\mathcal P}(\mu),
\end{equation*}
where the disjoint union runs over all the partitions of 
$M=0,1,\ldots, [L/2]$.
${\mathcal P}(\mu)$ is the set of states whose configuration is $\mu$. 
Each subset ${\mathcal P}(\mu)$ is invariant 
under any time evolution $T_l$, telling us that 
$\mu$ is a conserved quantity (\cite{KTT} Corollary 3.5).
Physical meaning of $\mu$ is the {\em soliton content}, 
namely, the list of the amplitudes of the 
solitons involved in $p$.
In particular $g$ is the number of solitons.

Unless otherwise stated, we shall consider 
those states whose configuration has the distinct parts as 
\begin{equation}\label{eq:mu}
\mu = \{i_1 < i _2 \cdots < i_g\}.
\end{equation}
Define the $g\times g$ symmetric integer matrix $A=(A_{i,j})_{i,j \in \mu}$ 
and the lattice $\Gamma$ by
\begin{equation}\label{eq:A}
A_{i,j} = \delta_{i,j}p_i + 2\min(i,j), \qquad \Gamma = A\Z^g \subset \Z^g.
\end{equation}
This matrix has arisen in the Bethe equation at $q=0$ (\ref{eq:sce})
known as the string centre equation \cite{KN}.
Under the condition $L \ge 2M$, $A$ is positive definite.

Let us proceed to the scattering data, i.e., the action-angle variables.
The action variable is the set $\mu$ itself.
The set of angle variables with prescribed $\mu$ 
is given by the quotient:
\begin{equation}\label{eq:udj}
{\mathcal J} = {\mathcal J}(\mu) = \Z^g/\Gamma.
\end{equation}
The one to be assigned with the state $p$ 
is found by the direct scattering map:
\begin{equation}\label{eq:ds}
\begin{split}
\Phi: \quad & {\mathcal P}(\mu)\;
\longrightarrow \;\Z \times {\mathcal P}_+
\longrightarrow \;\;\;\;
{\mathcal J}(\mu)\\
&\;\;p \;\quad\,\longmapsto \;(d, p_+) \,\longmapsto\;
({\bf J} + d{\bf h}_1)/\Gamma,
\end{split}
\end{equation}
where ${\bf h}_1 = (1,\ldots, 1) \in \Z^g$ as defined in (\ref{eq:tlj}).
${\bf J} = (J_i)_{i \in \mu} \in \Z^g$ is specified by the KKR bijection as in 
(\ref{eq:phi}), which we  
write as $\phi(p_+) = (\mu, {\bf J})$ or simply 
$\phi(p_+) = {\bf J}$.
Then ${\bf J} + d{\bf h}_1 = (J_i + d)_{i \in \mu}$.
$\Phi$ is well-defined \cite{KTT}. 
In particular,  the non-uniqueness of the decomposition 
$p \mapsto (d,p_+)$ is cancelled by taking 
${\rm mod } \;\Gamma$. 
For ${\bf I} \in \Z^g$, we denote its image in ${\mathcal J}$ 
by the same symbol ${\bf I}$.

For ${\bf I} \in {\mathcal J}$ we introduce the time evolution 
through 
\begin{equation}\label{eq:tlj}
T_l({\bf I}) = {\bf I} + {\bf h}_l,\quad 
{\bf h}_l = (\min(i,l))_{i\in \mu} \in \Z^g.
\end{equation}
Note that 
$L{\bf h}_1 = A{\bf h}_1 \in \Gamma$, therefore
$T_1^L({\bf I}) = {\bf I} \in {\mathcal J}$.

\begin{theorem}[\cite{KTT}, Theorems 3.11, 3.12]\label{th:ktt}
The map $\Phi$ is a bijection and the following commutative diagram is valid:
\begin{equation}\label{eq:cd}
\begin{CD}
{\mathcal P}(\mu) @>{\Phi}>> {\mathcal J}(\mu) \\
@V{T_l}VV @VV{T_l}V\\
{\mathcal P}(\mu) @>{\Phi}>> {\mathcal J}(\mu) 
\end{CD}
\end{equation}
Here $T_l$ on the left and the right are given by 
(\ref{eq:tl}) and (\ref{eq:tlj}), respectively.
\end{theorem}
The composition $\Phi^{-1}\circ T_l \circ \Phi$ 
yields the algorithmic solution of the initial value problem 
by the inverse scattering method \cite{GGKM,AS}. 
The nonlinear dynamics $T_l$ on ${\mathcal P}(\mu)$ 
becomes the straight motion on ${\mathcal J}(\mu)$ 
with the velocity ${\bf h}_l$. 
In this sense ${\mathcal J}(\mu)$ is an ultradiscrete 
analogue of the Jacobi variety.
Its cardinality is given by 
$\vert {\mathcal J}(\mu) \vert = \det A = Lp_{i_1}p_{i_2} \cdots p_{i_{g-1}}$
(\cite{KTT}, (4.6),(4.13) and (4.21)).
For $l \ge i_g$, one has ${\bf h}_l = {\bf h}_{i_g}$, hence  
$T_l(p) = T_{i_g}(p)$ by theorem \ref{th:ktt}.

In the limit $L \rightarrow \infty$, the quotient by $\Gamma$
in (\ref{eq:udj}) becomes void and the result provides the 
inverse scattering method for the box-ball system 
on the infinite lattice. 
The direct and the inverse scattering maps $\Phi^{\pm1}$ 
reduce to the KKR bijection $\phi^{\pm 1}$ itself. 

\begin{example}
For $p = 22121111222111$, let us derive 
\begin{equation}\label{eq:ans}
T^{1000}_2(p) = 11112221112212,\qquad
T^{1000}_3(p) = 12211122111122
\end{equation}
based on the inverse scattering scheme (\ref{eq:cd}).
(This $p$ is $T_2(p)$ in example \ref{ex:t23}.)
We have $p = T_1^2(p_+)$ with the highest state  
$p_+ = 12111122211122$.
The image of the KKR bijection of $\phi(p_+)$ and the
direct scattering transform $\Phi(p)$ are given by 

\begin{picture}(60,50)(-80,-40)
\put(-42,-19){$\phi(p_+) = $}
\put(0,0){\line(1,0){30}} \put(33,-9){1}
\put(0,-10){\line(1,0){30}}\put(23,-19){4}
\put(0,-20){\line(1,0){20}}\put(13,-29){0}
\put(0,-30){\line(1,0){10}}

\put(0,0){\line(0,-1){30}}
\put(10,0){\line(0,-1){30}}
\put(20,0){\line(0,-1){20}}
\put(30,0){\line(0,-1){10}}

\put(108,-19){$\Phi(p) = $}
\put(150,0){\line(1,0){30}} \put(183,-9){3}
\put(150,-10){\line(1,0){30}}\put(173,-19){6}
\put(150,-20){\line(1,0){20}}\put(163,-29){2}
\put(150,-30){\line(1,0){10}}

\put(150,0){\line(0,-1){30}}
\put(160,0){\line(0,-1){30}}
\put(170,0){\line(0,-1){20}}
\put(180,0){\line(0,-1){10}}

\end{picture}

\noindent
Thus $\mu=\{1,2,3\}, (p_1,p_2,p_3) = (8,4,2)$ 
and the matrix $A$ (\ref{eq:A})
reads
\begin{equation*}
A = \begin{pmatrix}
p_1+2 & 2 & 2 \\
2 & p_2+4 & 4 \\
2 & 4 & p_3+6
\end{pmatrix}
= \begin{pmatrix}
10 & 2 & 2 \\
2 & 8 & 4 \\
2 & 4 & 8
\end{pmatrix}.
\end{equation*}
According to (\ref{eq:cd}) and (\ref{eq:tlj}), the scattering data for the 
states $T^{1000}_{2,3}(p)$ are given by

\begin{picture}(60,50)(-80,-40)

\put(-68,-19){$T^{1000}_2(\Phi(p)) = $}
\put(0,0){\line(1,0){30}} \put(33,-9){2003}
\put(0,-10){\line(1,0){30}}\put(23,-19){2006}
\put(0,-20){\line(1,0){20}}\put(13,-29){1002}
\put(0,-30){\line(1,0){10}}

\put(0,0){\line(0,-1){30}}
\put(10,0){\line(0,-1){30}}
\put(20,0){\line(0,-1){20}}
\put(30,0){\line(0,-1){10}}

\put(82,-19){$T^{1000}_3(\Phi(p)) =$}
\put(150,0){\line(1,0){30}} \put(183,-9){3003}
\put(150,-10){\line(1,0){30}}\put(173,-19){2006}
\put(150,-20){\line(1,0){20}}\put(163,-29){1002}
\put(150,-30){\line(1,0){10}}

\put(150,0){\line(0,-1){30}}
\put(160,0){\line(0,-1){30}}
\put(170,0){\line(0,-1){20}}
\put(180,0){\line(0,-1){10}}

\end{picture}

\noindent
The angle variables appearing here are written as 
\begin{equation*}
\begin{pmatrix}
1002\\ 2006 \\ 2003
\end{pmatrix} 
=
\begin{pmatrix}
8\\ 4 \\ 1
\end{pmatrix} 
+ 0 {\bf h}_1 +
A
\begin{pmatrix}
35 \\ 161\\ 161
\end{pmatrix},
\qquad
\begin{pmatrix}
1002\\ 2006 \\ 3003
\end{pmatrix} 
=
\begin{pmatrix}
6\\ 0 \\ 1
\end{pmatrix} +
4 {\bf h}_1 + 
A
\begin{pmatrix}
17 \\ 81\\ 330
\end{pmatrix}.
\end{equation*}
The last terms involving $A$ can be dropped by 
${\rm mod }\, \Gamma$, whereas 
the first terms in the right hand sides give rise to 
the rigged configurations and the corresponding highest states:

\begin{picture}(60,50)(-120,-40)

\put(-100,-19)
{11112221112212$\;\;\overset{\;\,\phi^{-1}}{\longleftarrow}$}
\put(0,0){\line(1,0){30}} \put(33,-9){1}
\put(0,-10){\line(1,0){30}}\put(23,-19){4}
\put(0,-20){\line(1,0){20}}\put(13,-29){8}
\put(0,-30){\line(1,0){10}}

\put(0,0){\line(0,-1){30}}
\put(10,0){\line(0,-1){30}}
\put(20,0){\line(0,-1){20}}
\put(30,0){\line(0,-1){10}}

\put(50,-19)
{11221111221221$\;\;\overset{\;\,\phi^{-1}}{\longleftarrow}$}
\put(150,0){\line(1,0){30}} \put(183,-9){1}
\put(150,-10){\line(1,0){30}}\put(173,-19){0}
\put(150,-20){\line(1,0){20}}\put(163,-29){6}
\put(150,-30){\line(1,0){10}}

\put(150,0){\line(0,-1){30}}
\put(160,0){\line(0,-1){30}}
\put(170,0){\line(0,-1){20}}
\put(180,0){\line(0,-1){10}}

\end{picture}

\noindent
In view of $+0{\bf h}_1$ and $+4{\bf h}_1$, 
$T_2^{1000}(p)$ and $T_3^{1000}(p)$ are obtained by
taking the cyclic shifts $T^0_1$ and $T^4_1$ of these states
respectively, in agreement with (\ref{eq:ans}).
\end{example}  

\section{The explicit formula for the initial value problem}\label{sec:main}

First we present a piecewise linear formula for the KKR bijection. 
Let $(\mu,{\bf J})$ be a rigged configuration for a highest state 
in $B^{\otimes L}_1$.
To be concrete, we set 
\begin{equation*}
\phi^{-1}((\mu,{\bf J})) = 
(1-y(1),y(1))\otimes \cdots \otimes (1-y(L),y(L)) 
\in {\mathcal P}_+,
\end{equation*}
where $y(k)\in \{0,1\}$ is the `number of balls' in the 
$k$ th box from the left.
We parametrize the configuration 
$\mu=\{i_1,\ldots, i_g\}$ and 
the rigging ${\bf J}=(J_{i_1},\ldots, J_{i_g})$ 
as in (\ref{eq:phi}).
The following proposition \ref{pr:kkr} and 
lemma \ref{lem:waru} hold for the configurations 
such that $i_1 \le \cdots \le i_g$.

\begin{proposition}\label{pr:kkr}
The image of the KKR bijection is given by
\begin{align}
y(k) &= \tau_0(k)-\tau_0(k-1)-\tau_1(k)+\tau_1(k-1),
\label{eq:xk}\\
\tau_r(k) &= -\min_{{\bf n} \in \{0,1\}^g}
\{\sum_{i \in \mu}(J_i+r i -k)n_i 
+ \sum_{i,j \in \mu}\min(i,j)n_in_j\}
\quad (r=0,1),\label{eq:tau}
\end{align}
where ${\bf n} = (n_{i_1},\ldots, n_{i_g})$.
\end{proposition}
The proof will be given elsewhere for a more general case.
$\tau_r(k) \in \Z_{\ge 0}$ is the ultradiscrete tau function   
mentioned in section \ref{sec:intro}.
We remark that there is {\em no} dependence on $L$ 
in (\ref{eq:tau}) except in 
the upper bound $p_i$ (\ref{eq:pi}) of the rigging
$J_i \le p_i$.
For $k<1$ or $k >L$, (\ref{eq:xk}) gives $y(k)=0$.
As it turns out, after theorem \ref{th:main},
proposition \ref{pr:kkr} essentially provides the solution
of the initial value problem of the box-ball system 
on the infinite lattice $k \in \Z$.

\begin{lemma}[\cite{KTT}, Lemma C.1]\label{lem:waru}
Let $q \in B^{\otimes K}_1$ and 
$r \in B^{\otimes L}_1$ be the highest states associated with
the rigged configurations 
$\phi(q) = (\lambda, {\bf I})$ and $\phi(r)= (\mu, {\bf J})$.
Then the rigged configuration of the 
highest state $q\otimes r \in B^{\otimes K+ L}_1$ is 
$\phi(q \otimes r) = (\lambda \cup \mu, {\bf I} \cup {\bf J}')$, where 
${\bf J}'=(J'_j)_{j\in \mu}$ is given by
\begin{equation*}
J'_j = J_j + p_j,\qquad 
p_j = K-2\sum_{k\in \lambda}\min(j,k).
\end{equation*}
\end{lemma}
The shift $p_j$ here is nothing but 
the vacancy number in the rigged configuration $\phi(q)$. 
The notation 
$(\lambda \cup \mu, {\bf I} \cup {\bf J}')$ means the union regarding 
$(\lambda, {\bf I})$ and $(\mu, {\bf J}')$ 
as multi-sets of parts (rows in Young diagrams) 
assigned with rigging. 
For example,

\begin{picture}(100,60)(-50,-45)
\put(-40,-10){$(\lambda, {\bf I}) = $}
\put(0,0){\line(1,0){20}}\put(23,-9){$a$}
\put(0,-10){\line(1,0){20}}\put(13,-19){$b$}
\put(0,-20){\line(1,0){10}}
\put(0,0){\line(0,-1){20}}
\put(10,0){\line(0,-1){20}}
\put(20,0){\line(0,-1){10}}

\put(55,-10){$(\mu, {\bf J}') = $}
\put(100,0){\line(1,0){30}}\put(133,-9){$c$}
\put(100,-10){\line(1,0){30}}\put(113,-19){$d$}
\put(100,-20){\line(1,0){10}}
\put(100,0){\line(0,-1){20}}
\put(110,0){\line(0,-1){20}}
\put(120,0){\line(0,-1){10}}
\put(130,0){\line(0,-1){10}}

\put(160,-10){$(\lambda \cup \mu, {\bf I} \cup {\bf J}') = $}
\put(240,0){\line(1,0){30}}\put(273,-9){$c$}
\put(240,-10){\line(1,0){30}}\put(263,-19){$a$}
\put(240,-20){\line(1,0){20}}\put(253,-29){$d$}
\put(240,-30){\line(1,0){10}}\put(253,-39){$b$}
\put(240,-40){\line(1,0){10}}

\put(240,0){\line(0,-1){40}}
\put(250,0){\line(0,-1){40}}
\put(260,0){\line(0,-1){20}}
\put(270,0){\line(0,-1){10}}

\end{picture}  

\noindent
where, as usual, the ordering of the rigging 
$d$ and $b$ 
within a block of equal length rows 
does not matter.
In what follows, we employ the convention of always arranging
the rigging to weakly increase upward within such blocks.

Given a state $p \in {\mathcal P}$, 
take a highest state $p_+ \in {\mathcal P}_+$ 
and $0 \le d <L$ such that $p=T^d_1(p_+)$.
Let $\phi(p_+) = (\mu,{\bf J})$ be the rigged configuration for $p_+$,
which we parametrize as 
$\mu=\{i_1,\ldots, i_g\}$ and ${\bf J}=(J_{i_1},\ldots, J_{i_g})$.
Here we assume $i_1 < \cdots < i_g$ in accordance with the 
assumption (\ref{eq:mu}).
We form a large highest state 
$p_+^{\otimes N} = p_+ \otimes \cdots \otimes p_+ 
\in B^{\otimes NL}_1$.
By lemma \ref{lem:waru}, its rigged configuration
$(\mu^N, {\bf J}^N) := \phi(p_+^{\otimes N})$ is given by
\begin{align*}
\mu^N & = \{i_{1,1}, \ldots, i_{1,N},
i_{2,1}, \ldots, i_{2,N}, \ldots, 
i_{g,1}, \ldots, i_{g,N}\},\\
{\bf J}^N & = (J_{i_1,1},\ldots,J_{i_1,N},
J_{i_2,1},\ldots,J_{i_2,N}, \ldots, 
J_{i_g,1},\ldots, J_{i_g,N}),\\
i_{s,\alpha} &= i_s,\qquad 
J_{i_s,\alpha} = J_{i_s} + (\alpha-1)p_{i_s}\quad
(1 \le \alpha \le N),
\end{align*}
where $p_i = L-2\sum_{j\in \mu}\min(i,j)$ is the vacancy number 
for $p_+$.
We apply proposition \ref{pr:kkr} to $(\mu^N,{\bf J}^N)$.
{}From (\ref{eq:tau}) 
the corresponding ultradiscrete tau function $\tau_r(k)$ reads
\begin{equation}\label{eq:chu}
-\min_{{\bf n} \in \{0,1\}^{Ng}}
\left\{\sum_{i \in \mu}\sum_{1\le \alpha \le N}
(J_{i,\alpha}+ri-k)n_{i,\alpha}
+ \sum_{i,j \in \mu}\sum_{1 \le \alpha, \beta \le N}
\min(i,j)n_{i,\alpha}n_{j,\beta} \right\},
\end{equation}
where 
${\bf n} = 
(n_{i_1,1},\ldots, n_{i_1,N},\ldots, n_{i_g,1},\ldots, n_{i_g,N})$.
Since $J_{i,1}\le J_{i,2} \le \cdots \le J_{i,N}$ for each $i \in \mu$, 
the minimum here can be restricted to those 
${\bf n}$ having the form
\begin{equation*}
n_{i,1}=n_{i,2}=\cdots = n_{i,m_i}=1,\quad
n_{i,m_i+1} = n_{i,m_i+2} = \cdots = n_{i,N} = 0
\end{equation*}
for some $0 \le m_i \le N$.
Then the sums over $\alpha$ and $\beta$ in (\ref{eq:chu})
can be taken, leading to 
\begin{equation}\label{eq:qf}
\begin{split}
&\sum_{i\in \mu}
\Bigl(m_iJ_i+\frac{m_i(m_i-1)}{2}p_i + 
m_iri - m_ik\Bigr) + \sum_{i,j \in \mu}\min(i,j)m_im_j\\
&=
{}^t{\bf m}\bigl({\bf J} - \frac{\bf p}{2}+ r{\bf h}_\infty
-k{\bf h}_1\bigr) + \frac{1}{2}{}^t{\bf m}A{\bf m},
\end{split}
\end{equation}
where $A=(A_{i,j})$ is defined in (\ref{eq:A}).
We have set 
${\bf m} = (m_i)_{i \in \mu}$, 
${\bf p} = (p_i)_{i \in \mu}$ and used 
the vector notation ${\bf J}, {\bf h}_1, {\bf h}_\infty$
around (\ref{eq:ds}) and (\ref{eq:tlj}). 
For instance 
${\bf h}_\infty = {\bf h}_{i_g}$ and 
(\ref{eq:pi}) is rephrased as
\begin{equation}\label{eq:pp}
{\bf p} = L{\bf h}_1 - 2 \sum_{j\in \mu}{\bf h}_j.
\end{equation}
By taking $N$ to be even and 
shifting ${\bf m}$ to ${\bf m}+\frac{N}{2}{\bf h}_1$,
(\ref{eq:qf}) is rewritten as
${}^t{\bf m}({\bf J} - \frac{\bf p}{2}+ r{\bf h}_\infty
-(k-\frac{NL}{2}){\bf h}_1) + \frac{1}{2}{}^t{\bf m}A{\bf m} + X$,
where 
$X = \frac{N}{2}{}^t{\bf h}_1
({\bf J} - \frac{\bf p}{2}+ r{\bf h}_\infty
-(k-\frac{NL}{4}){\bf h}_1)$.
This $X$ can be put outside $\min$, after which 
its dependence on $r,k$ is cancelled 
in the difference (\ref{eq:xk}).
Therefore we find that 
$p_+^{\otimes N} = 
(1-y(1),y(1))\otimes \cdots \otimes (1-y(NL),y(NL))$ is 
given by (\ref{eq:xk}) with $\tau_r(k)$ replaced by
\begin{equation}\label{eq:bakka}
\tau_r(k) = -\min_{\bf m}\{
{}^t{\bf m}({\bf J} - \frac{\bf p}{2}+ r{\bf h}_\infty
-(k-\frac{NL}{2}){\bf h}_1) + \frac{1}{2}{}^t{\bf m}A{\bf m}\},
\end{equation}
where $\min$ is taken over those 
${\bf m} =(m_i)_{i \in \mu} \in \Z^g$ 
such that $-N/2 \le m_i \le N/2$.

{}From the relation $p = T^d_1(p_+)$, 
the state 
$p=(1-x(1),x(1))\otimes \cdots \otimes (1-x(L),x(L))$ 
is obtained {}from $p_+^{\otimes N}$ 
by picking up the length $L$ segment 
corresponding to 
$y(wL-d+1), \ldots, y((w+1)L-d)$
for any  $1 \le w \le N-1$.
Thus in (\ref{eq:bakka}) 
we replace $k$ by $k + wL-d$ 
with the choice $w=\frac{N}{2}$ to get 
$\tau_r(k) = -\min_{\bf m}\{c_L({\bf m})\}$ with
\begin{equation}\label{eq:cl}
c_L({\bf m}) = {}^t{\bf m}
\bigl({\bf I}-\frac{\bf p}{2}-k{\bf h}_1+r{\bf h}_\infty\bigr)
+ \frac{1}{2}{}^t{\bf m}A{\bf m}.
\end{equation}
Here we have let 
${\bf I}= {\bf J} + d{\bf h}_1$ 
denote the angle variable $\Phi(p)$ for $p$. 
See (\ref{eq:ds}).
The resulting formula for $x(k)$
gives the state $p$ corresponding to its 
action-angle variable
$(\mu, \,{\bf I})$
{\em as long as} $0 \le d \le L-1$, 
$1 \le k \le L$ and $0\le J_i \le p_i$
since we have started from the rigged configuration.
These constraints are removed by taking the 
limit $N \rightarrow \infty$, where the minimum 
extends over ${\bf m} \in \Z^g$; therefore one has
\begin{equation*}
\tau_r(k) = \Theta\bigl(
{\bf I}-\frac{\bf p}{2}-k{\bf h}_1+r{\bf h}_\infty\bigr).
\end{equation*}
By virtue of the quasi-periodicity of the 
ultradiscrete Riemann theta function (\ref{eq:qp}),
the difference 
\begin{equation}\label{eq:y}
\begin{split}
x(k) &= \Theta\bigl(
{\bf I}-\frac{\bf p}{2}-k{\bf h}_1\bigr) - 
\Theta\bigl(
{\bf I}-\frac{\bf p}{2}-(k\!-\!1){\bf h}_1\bigr)\\
&-\Theta\bigl(
{\bf I}-\frac{\bf p}{2}-k{\bf h}_1+{\bf h}_\infty\bigr)+
\Theta\bigl(
{\bf I}-\frac{\bf p}{2}-(k\!-\!1){\bf h}_1+{\bf h}_\infty\bigr)
\end{split}
\end{equation}
gains the invariance under 
$k \rightarrow k+L$ 
and ${\bf I} \rightarrow {\bf I} + {\bf v}$ for any 
${\bf v} \in \Gamma = A\Z^g$.
(Note that $L{\bf h}_1 = A{\bf h}_1 \in \Gamma$.)
Namely, (\ref{eq:y}) makes sense for 
$k \in \Z_L$ and ${\bf I} \in {\mathcal J}=\Z^g/\Gamma$.

To summarize, we have proved 
\begin{theorem}\label{th:main}
For any state $p \in {\mathcal P}$ of the 
periodic box-ball system,
let $(\mu,\, {\bf I}) = \Phi(p)$ be the action-angle variable.
Fix ${\bf p}=(p_i)_{i \in \mu}$ by (\ref{eq:pp}) and 
the matrix $A$ by (\ref{eq:A}).
Then the state $p$ is expressed as 
$p=(1-x(1),x(1))\otimes \cdots \otimes (1-x(L),x(L))$
with $x(k)\in \{0,1\}$ given by (\ref{eq:y}).
Due to theorem \ref{th:ktt}, 
this solves the initial value problem 
in that any time evolution 
$T_{l_1}^{\gamma_1}\cdots T_{l_t}^{\gamma_t}(p)$ 
is obtained by replacing ${\bf I}$ in (\ref{eq:y}) with 
${\bf I} + \gamma_1{\bf h}_{l_1}+ \cdots + \gamma_t{\bf h}_{l_t}
\; (\gamma_i \in \Z)$.
\end{theorem}

The quadratic form (\ref{eq:cl}) is decomposed as
$c_L({\bf m}) = L\sum_{i=1}^gm_i(m_i-1)/2 + c({\bf m})$, where
$c({\bf m})$ is independent of the system size $L$.
In the limit $L \rightarrow \infty$,
the minimum is restricted to ${\bf m} \in \{0,1\}^g$ and 
$\Theta$ degenerates into the ultradiscrete tau function
as in the scheme (\ref{eq:4}).
If ${\bf I}$ is chosen to be a rigged configuration,
the formula (\ref{eq:y}) under such a reduction 
still describes the image of the KKR bijection although 
the function $-\min_{{\bf m} \in \{0,1\}^g}\{c({\bf m})\}$ 
takes slightly different form from (\ref{eq:tau}).
The result provides the solution of the initial value 
problem of the box-ball system on the infinite lattice.

In Figure 1, we plot the following function
on the $(k,t)$ (space-time) plane:
\begin{equation}\label{eq:ty}
u(k,t) = \frac{
\vartheta\bigl(
T^t_{\infty}({\bf I})-\frac{\bf p}{2}-k{\bf h}_1\bigr)
\vartheta\bigl(
T^t_{\infty}({\bf I})-\frac{\bf p}{2}-(k\!-\!1){\bf h}_1+{\bf h}_\infty\bigr)}
{\vartheta\bigl(
T^t_{\infty}({\bf I})-\frac{\bf p}{2}-(k\!-\!1){\bf h}_1\bigr)
\vartheta\bigl(
T^t_{\infty}({\bf I})-\frac{\bf p}{2}-k{\bf h}_1+{\bf h}_\infty\bigr)},
\end{equation}
where $T^t_\infty({\bf I}) = {\bf I} + t{\bf h}_\infty$ by 
(\ref{eq:tlj}) and 
$\vartheta({\bf z}) = \sum_{{\bf n} \in \Z^g}
\exp\Bigl(-({}^t{\bf n}A{\bf n}/2+{}^t{\bf n}{\bf z})/\epsilon
\Bigr)$
is the Riemann theta function.
In view of the scheme (\ref{eq:4}), one has 
$\lim_{\epsilon\rightarrow +0}\epsilon\log u(k,0) = x(k)$.
Thus $u(k,t)$ gives a softening 
of the envelop of ultradiscrete solitons in 
the periodic box-ball system at $\epsilon = 0$ under the 
time evolution $T_\infty$.
The selected parameters are
\begin{equation*}
L=170, \;\mu = \{2,6\},
{\bf I} = \begin{pmatrix}0 \\ 0 \end{pmatrix},\;
{\bf p} = \begin{pmatrix}p_2 \\ p_6 \end{pmatrix}
= \begin{pmatrix} 162 \\ 154 \end{pmatrix},\;
A=\begin{pmatrix}166 & 4 \\ 4 & 166 \end{pmatrix},\;
\epsilon= 7.
\end{equation*}
For the periodic box-ball system described by (\ref{eq:y}),
this data corresponds to
$p = 1122111111222222\otimes 1^{\otimes 154}$, which is a 
two soliton state with amplitudes 2 and 6.
At $t=70$, it becomes 
$T^{70}_\infty(p) =  1^{\otimes 94}\otimes 222222 \otimes 
1^{\otimes 38} \otimes 22 \otimes 1^{\otimes 30}$.

\begin{figure}[h]

\includegraphics{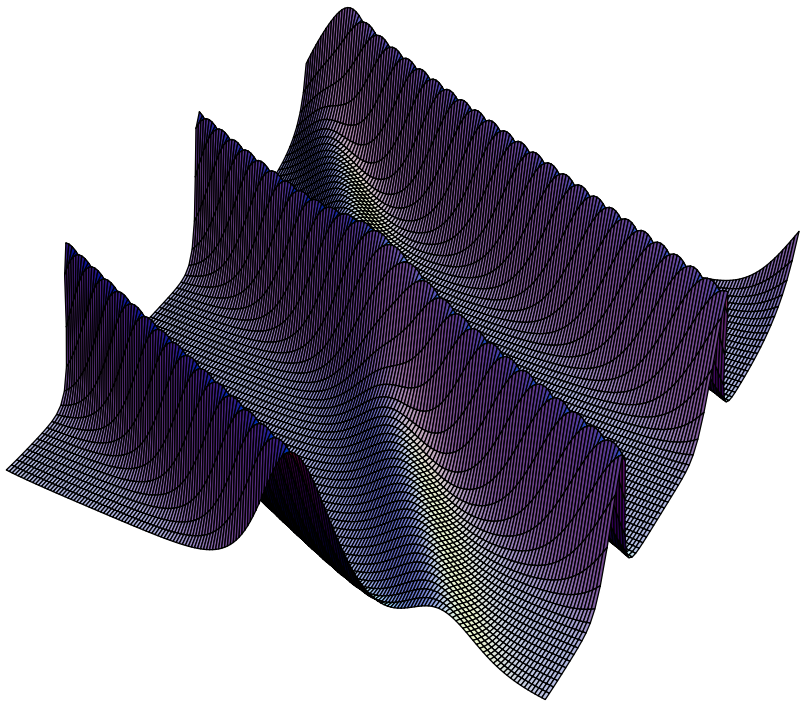}

\caption{The envelope of the function 
$u(k,t)$ (\ref{eq:ty}) 
for $1 \le k \le 170,\, 0 \le t \le 70$.
The top and right corners 
correspond to $(k,t)=(0,0), (170,0)$, respectively.
It is periodic in the $k$-direction.}
\end{figure}

\section{Discussion}\label{sec:discussion} 
Theorem \ref{th:main} enables one to construct the joint eigenvectors 
of $T_1, T_2, \ldots$ in $({\mathbb C}^2)^{\otimes L}$.
The result may be regarded as an 
explicit formula for $q=0$ Bethe vectors 
in terms of the ultradiscrete Riemann theta function.
We continue assuming that $\mu$ consists of distinct parts as in 
(\ref{eq:mu}).

The Bethe equation for 
the periodic XXZ chain on 
$({\mathbb C}^2)^{\otimes L}$ 
associated with $U_q(\widehat{sl}_2)$ 
becomes linear at $q=0$ under the string hypothesis.
The result is known as the string centre equation \cite{KN}:
\begin{equation}\label{eq:sce}
A{\bf u} \equiv -\frac{\bf p}{2} \mod \Z^g,
\end{equation}
where ${\bf u}=(u_{i_1},\ldots, u_{i_g})$ 
with $u_i$ being the centre of the length $i$ string.
We call ${\bf u}$ the Bethe root.
In this normalization, the Bethe wave function is a 
rational function of $\exp(2\pi\sqrt{-1}u_i)$; hence
${\bf u}$ lives in $(\R/\Z)^g$. 
Thus there is one to one correspondence between the 
Bethe root ${\bf u}$ and 
the angle variable ${\bf J} \in {\mathcal J} = \Z^g/A\Z^g$ 
via the relation \cite{KTT}
\begin{equation*}
A{\bf u} = {\bf J} - \frac{\bf p}{2}.
\end{equation*}
The time evolution $T_l$ of ${\bf J}$ (\ref{eq:tlj})
induces that of the Bethe roots, which is again 
a straight motion $T_l({\bf u}) = {\bf u} + A^{-1}{\bf h}_l$
in $(\R/\Z)^g$. 

At first sight, this appears contradictory, because 
$T_1, T_2, \ldots$ are fusion transfer matrices at $q=0$, 
which should leave the $q=0$ Bethe vectors invariant 
up to an overall scalar as well as the relevant Bethe roots.
The answer to this puzzle is that 
the state $p \in B^{\otimes L}_1$ that we are 
associating to ${\bf u}$ or ${\bf J}$ 
by $\Phi(p) = (\mu,{\bf J})$ is a {\em monomial} 
in $({\mathbb C}^2)^{\otimes L}$, which is {\em not} 
a Bethe vector at $q=0$ in general.

It is easy to remedy this.
In fact, for each Bethe root ${\bf u}$ or equivalently  
${\bf J}=A{\bf u} + \frac{\bf p}{2} \in {\mathcal J}$,
one can construct a vector 
$\vert {\bf J} \rangle  \in ({\mathbb C}^2)^{\otimes L}$ 
that possesses every aspect as a $q=0$ Bethe vector as follows:
\begin{align}
\vert {\bf J}\rangle &= \sum_{{\bf I} \in {\mathcal J}}
c_{{\bf I}, {\bf J}}\,p({\bf I}),\label{eq:jp}\\
c_{{\bf I}, {\bf J}} &= \exp\left(-2\pi\sqrt{-1}\;
{}^t{\bf I}\Bigl(A^{-1}({\bf J}-\frac{\bf p}{2})
+\frac{{\bf h}_1}{2}\Bigr)\right), \nonumber\\
p({\bf I}) & =
\begin{pmatrix}1-x(1)\\ x(1) \end{pmatrix}
\otimes \cdots \otimes 
\begin{pmatrix}1-x(L) \\ x(L) \end{pmatrix}
\in {\mathcal P}(\mu) \subseteq ({\mathbb C}^2)^{\otimes L}, 
\nonumber
\end{align}
where $x(k)\in \{0,1\}$ is specified by (\ref{eq:y}).
We embed $B^{\otimes L}_1$ into $({\mathbb C}^2)^{\otimes L}$
naturally and extend $T_l$ to the latter by ${\mathbb C}$-linearity.
The vector $p({\bf I})$ here is nothing but the state of the 
periodic box-ball system appearing in theorem \ref{th:main}.
It follows that
$T_l(p({\bf I})) = p({\bf I}+{\bf h}_l)$.
Thus from ${\mathcal J} + {\bf h}_l = {\mathcal J}$, 
it is elementary to check
\begin{align*}
T_l\vert {\bf J}\rangle &= \Lambda_l({\bf J})\vert {\bf J}\rangle,\\
\Lambda_l({\bf J}) & = c_{-{\bf h}_l, {\bf J}}
= \exp\left(2\pi\sqrt{-1}\;
{}^t{\bf h}_l\Bigl({\bf u}+\frac{{\bf h}_1}{2}\Bigr)\right).
\end{align*}
The quantity $\Lambda_l({\bf J})$ here exactly coincides 
with the $q=0$ Bethe eigenvalue given in equation (4.28) of \cite{KTT}.
Note further that the transition relation (\ref{eq:jp}) is inverted as
\begin{equation*}
p({\bf I}) = \frac{1}{\vert {\mathcal J} \vert}
\sum_{{\bf J} \in {\mathcal J}} {\bar c}_{{\bf I}, {\bf J}} 
\vert {\bf J}\rangle,
\end{equation*}
where ${\bar c}_{{\bf I}, {\bf J}}$ denotes the complex conjugate 
of $c_{{\bf I}, {\bf J}}$. 
It follows that the space of the $q=0$ Bethe vectors 
$\vert {\bf J} \rangle$  
coincides with the space of periodic box-ball states $p$ for each 
prescribed soliton content $\mu$, namely,
\begin{equation*}
\bigoplus_{{\bf J}  \in {\mathcal J}(\mu)}{\mathbb C}
\vert {\bf J} \rangle
= \bigoplus_{p \in {\mathcal P}(\mu)} {\mathbb C} \,p.
\end{equation*}

Thus we conclude that the approach here bypasses 
the formidable task of computing the $q\rightarrow 0$ limit of the
Bethe vectors in general, but leads to the joint eigenvectors   
$\vert {\bf J} \rangle$ of $\{T_l\}$.
They form a basis of the space having the 
prescribed soliton content and possess 
the spectrum $\Lambda_l({\bf J})$ anticipated 
from the Bethe ansatz at $q=0$. 
Moreover $\vert {\bf J} \rangle$ is 
parametrized explicitly 
in terms of the ultradiscrete Riemann theta function.

\vspace{0.3cm}\noindent
{\bf Acknowledgments} \hspace{0.1cm}
The authors thank Tomoki Nakanishi, Masato Okado, Mark Shimozono, 
Taichiro Takagi, Akira Takenouchi 
and Yasuhiko Yamada for discussion on 
related topics.
RS is a research fellow of the 
Japan Society for the Promotion of Science.
He thanks Miki Wadati for continuous encouragement.

\vspace{5mm}
\begin{flushleft}
Atsuo Kuniba:\\
Institute of Physics, Graduate School of Arts and Sciences,
University of Tokyo,
Komaba, Tokyo 153-8902, Japan\\
\texttt{atsuo@gokutan.c.u-tokyo.ac.jp}\vspace{3mm}\\
Reiho Sakamoto:\\
Department of Physics, Graduate School of Science, 
University of Tokyo, Hongo, Tokyo 113-0033, Japan\\
\texttt{reiho@monet.phys.s.u-tokyo.ac.jp}
\end{flushleft}

\begin{thebibliography}{A}


\bibitem{AS}
M. J. Ablowitz and H. Segur, 
Solitons and the inverse scattering transform,
SIAM Studies in Appl. Math. 4. Philadelphia Pa. (1981).

\bibitem{Ba}
R.~J.~Baxter,
Exactly solved models in statistical mechanics, 
Academic Press, London (1982).

\bibitem{Be}
H.\ A.\ Bethe,
Zur Theorie der Metalle, I. Eigenwerte und
Eigenfunktionen der linearen Atomkette,
Z.\ Physik {\bf 71} (1931) 205--231.

\bibitem{DT}
E. Date and S. Tanaka, 
Periodic multi-soliton solutions of Korteweg-de Vries equation 
and Toda lattice,
Prog. Theoret. Phys. Suppl.  {\bf 59}  (1976) 107--125.

\bibitem{DMN}
B. A. Dubrovin,  V. B. Matveev and S. P. Novikov, 
Nonlinear equations of Korteweg-de Vries type, 
finite-band linear operators and Abelian varieties
Russian Math. Surveys {\bf 31} (1976) 59--146.

\bibitem{GGKM}
C. S. Gardner, J. M. Greene, M. D. Kruskal and R. M. Miura,
Method for solving the Korteweg-de Vries equation,
Phys. Rev. Lett. {\bf 19} (1967) 1095--1097.

\bibitem{JM}
M. Jimbo and T. Miwa,
Solitons and infinite dimensional Lie algebras,
Publ. RIMS. Kyoto Univ. {\bf 19} (1983) 943--1001.

\bibitem{KR}
A. N. Kirillov and N. Yu. Reshetikhin, 
The Bethe ansatz and the combinatorics of Young tableaux. 
J. Soviet Math. {\bf 41} (1988) 925--955. 

\bibitem{KN}
A. Kuniba and T. Nakanishi, 
The Bethe equation at $q=0$, the M\"obius inversion formula, 
and weight multiplicities: I. 
The $sl(2)$ case, Prog. in Math. {\bf 191} (2000) 185--216.

\bibitem{KTT} 
A. Kuniba, T. Takagi and A. Takenouchi,
Bethe ansatz and inverse scattering transform 
in a periodic box-ball system, 
Nucl. Phys. B [PM] (2006) 354--397.

\bibitem{M} 
D. Mumford,
Tata Lectures on Theta II, Birkh\"auser, Boston (1984).

\bibitem{YYT}
D. Yoshihara, F. Yura and T. Tokihiro,
Fundamental cycle of a periodic box-ball system,
J. Phys. A: Math. Gen. {\bf 36}  (2003) 99--121.

\end{thebibliography}
\end{document}